\documentclass[12pt,a4paper,oneside,reqno,notitlepage]{amsart}
\usepackage{amscd,amssymb,hyperref}
\usepackage[PostScript=dvips]{diagrams}
\usepackage{graphicx}
\usepackage{graphics}
\usepackage{epsfig}
\usepackage{picinpar}
\usepackage{wrapfig}
\usepackage{amsmath}
\usepackage{color}
\usepackage{pstricks,pst-node,pst-text,pst-3d}
\usepackage{pst-plot}
\usepackage{verbatim}

\newtheorem{thm}{Theorem}[section]
\newtheorem{lem}[thm]{Lemma}

\newtheorem{prop}[thm]{Proposition}

\def\square{\vbox{
      \hrule height 0.4pt
      \hbox{\vrule width 0.4pt height 5.5pt \kern 5.5pt \vrule width 0.4pt}
      \hrule height 0.4pt}}

\def\Ker{\mathrm{K er}}
\def\ch\mathrm{c h}

\def\CP{\mathbb{C}\mathrm{P}}

\newcommand{\Z}{\mathbb{Z}}

\newcommand{\VAP}{\mathrm{VAP}}
\newcommand{\Brun}{\mathrm{Brun}}

\let\la=\langle
\let\ra=\rangle

\numberwithin{equation}{section}

\newcommand{\auths}[1]{\textrm{#1},}
\newcommand{\artTitle}[1]{\textsl{#1},}
\newcommand{\jTitle}[1]{\textrm{#1}}
\newcommand{\Vol}[1]{\textbf{#1}}
\newcommand{\Year}[1]{\textrm{(#1)}}
\newcommand{\Pages}[1]{\textrm{#1}}

\begin{document}

\title[homotopy aspects of virtual braids]{On homotopy aspects and cablings of virtual pure braid groups}

\author[V. Bardakov]{Valeriy G. Bardakov}
\address{Sobolev Institute of Mathematics, Novosibirsk 630090, Russia,}
\address{Novosibirsk State University, Novosibirsk 630090, Russia,}
\address{Novosibirsk State Agrarian University, Dobrolyubova street, 160,  Novosibirsk 630039, Russia,}
\email{bardakov@math.nsc.ru}

\author{Roman Mikhailov}
\address{Laboratory of Modern Algebra and Applications, St. Petersburg State University, 14th Line, 29b,
Saint Petersburg, 199178 Russia and St. Petersburg Department of
Steklov Mathematical Institute} \email{rmikhailov@mail.ru}

\author{Jie Wu }
\address{Department of Mathematics, National University of Singapore, 10 Lower Kent Ridge Road, Singapore 119076} \email{matwuj@nus.edu.sg}
\urladdr{www.math.nus.edu.sg/\~{}matwujie}
\thanks{The main result is supported by the Russian Science Foundation grant N 16-11-10073.}

\begin{abstract}
This is an announcement on the main results in~\cite{BMW1,BMW2}. By exploring simplicial structure of pure virtual braid groups, we give new connections between the homotopy groups of the $3$-sphere and the virtual braid groups that lead applications to the theory of Brunnian virtual braids~\cite{BMW1}. In~\cite{BMW2}, we prove that the group structure of $VP_n$ with $n\geq 5$ is determined by $VP_3$, $VP_4$ and virtual cablings given by iterated degeneracy operations on the generators and defining relations.

\end{abstract}

\maketitle

\section{Introduction} \label{int}
This is an announcement on the main results in~\cite{BMW1,BMW2} for exploring the homotopy aspects and the cablings of pure virtual braid groups with applications to the theory of Brunnian virtual braids. We outline the main statements in this article. The detailed exploration will be given in~\cite{BMW1,BMW2}.

\section{The homotopy aspects of virtual braid groups}\label{section2}
The exploration on simplicial structure on braids, links and mapping classes have been to lead some fundamental connections with homotopy groups~\cite{BCWW,BHW2,LLW,MW}. These connections enrich the meanings of homotopy groups in the content of the objects from low dimensional topology. In particular, it gives a description of elements in homotopy groups of the $2$-sphere in terms of Brunnian braids~\cite{BCWW} with a generalization on higher dimensional spheres in~\cite{MW}. The purpose of this article is to explore the simplicial structure on virtual braids. The basic ideas follow the approach introduced by F. Cohen and the third author~\cite{CW} on the classical braids briefly reviewed as follows.

The notion of simplicial set is a generalization of simplicial complex. A precise definition of simplicial set is as follows. A sequence of sets $X_* = \{ X_n \}_{n \geq 0}$  is called a
{\it simplicial set} if there are face maps:
$$
d_i : X_n \longrightarrow X_{n-1} ~\mbox{for}~0 \leq i \leq n
$$
and degeneracy  maps
$$
s_i : X_n \longrightarrow X_{n+1} ~\mbox{for}~0 \leq i \leq n,
$$
that are  satisfy the following simplicial identities:
\begin{enumerate}
\item $d_i d_j = d_{j-1} d_i$ if $i < j$,
\item $s_i s_j = s_{j+1} s_i$ if $i \leq j$,
\item $d_i s_j = s_{j-1} d_i$ if $i < j$,
\item $d_j s_j = id = d_{j+1} s_j$,
\item $d_i s_j = s_{j} d_{i-1}$ if $i > j+1$.
\end{enumerate}
Here $X_n$ can be geometrically viewed as the set of $n$-simplices including all possible degenerate simplices. For instance, a vertex $v$ can be viewed as a \textit{degenerate $1$-simplex} in the form of $[vv]$. An advantage of the notion of simplicial set is that one can do homotopy theory using this notion. A \textit{pointed simplicial set} is a simplicial set $X_*$ with a fixed choice of basepoint. In other words, we have a fixed basepoint $\ast\in X_0$ that creates one and only one degenerate $n$-simplex in $X_n$ by applying iterated degeneracy operations on it. A \textit{simplicial group} is a group object in the category of simplicial sets, namely a simplicial group is a simplicial set $X_*$ such that each $X_n$ is a group and all face and degeneracy operations are group homomorphisms. Milnor's $F[K]$ construction~\cite{Milnor} is the adjoint functor to the forgetful functor from the category of simplicial groups to the category of pointed simplicial sets. For a given pointed simplicial set $X_*$, the group $F[X_*]_n$ is the free group generated by $X_n$, the set of $n$-simplices, modulo the single relation that the basepoint equals to $1$. The geometric realization of $F[X_*]$ is homotopy equivalent to $\Omega \Sigma |X_*|$, the loop suspension of the geometric realization of $X_*$ ((Note. In Milnor's
paper~\cite{Milnor}, $X_*$ is required to be a reduced simplicial
set. This result actually holds for any pointed simplicial set by
a more general result~\cite[Theorem 4.9]{Wu0}.))

Now we review the construction in~\cite{CW}. Let $P_{n+1}$ be the group of the $(n+1)$-strand Artin pure braids with their strands labeled from $0$ to $n$. Consider the sequence of Artin pure braid groups $AP_*=\{AP_n\}_{n\geq0}$ with $AP_n=P_{n+1}$. The simplicial structure on $AP_*$ is given in such a way that the face operation $d_i\colon AP_n\to AP_{n-1}$ is given by deleting $i$-th strand, and the degeneracy operation $s_i\colon AP_n\to AP_{n+1}$ is given by doubling $i$-th strand. (Note. In this system, the $0$-th strand refers to the first strand.) Observe that $AP_0=P_1$ is the trivial group, and $AP_1=P_2=\mathbb{Z}$ generated by $A_{1,2}$. Consider $A_{1,2}\in AP_1$ as a $1$-simplex in the simplicial group $AP_*$. Since $d_0A_{1,2}=d_1A_{1,2}=1$, there is a unique simplicial map $f_{A_{1,2}}\colon S^1\longrightarrow AP_*$ such that $f_*((0,1))=A_{1,2}$. Here $S^1$ is the simplicial $1$-sphere with the non-degenerate $1$-simplex $(0,1)$. By using the universal property of Milnor's $F[K]$-construction, the simplicial map $f_{A_{1,2}}$ extends uniquely to a simplicial homomorphism
$$
\Theta\colon F[S^1]\longrightarrow AP_*
$$
with the image of $\Theta$ is the smallest simplicial subgroup of $AP_*$ containing $A_{1,2}$.
It was proved in ~\cite{CW} that $\Theta$ is injective, and so the simplicial group $F[S^1]$ embeds into $AP_*$ via $\Theta$.  Since the geometric realization of $F[S^1]$ is homotopy equivalent to the loop space $\Omega S^2$, the above representation describes $\pi_n(S^2)$ as a quotient of a subgroup of $P_n$.

Let $VP_{n+1}$ denote the virtual pure braid group on $n+1$ strands labeled from $0$ to $n$, where $VP_1=\{1\}$. Let $\VAP_n=VP_{n+1}$ for $n\geq0$. Along the methodology of  deleting-cabling strands as in~\cite{CW}, we obtain a simplicial group
$$
\VAP_* :\ \ \ \ldots\ \begin{matrix}\longrightarrow\\[-3.5mm] \ldots\\[-2.5mm]\longrightarrow\\[-3.5mm]
\longleftarrow\\[-3.5mm]\ldots\\[-2.5mm]\longleftarrow \end{matrix}\ VP_4 \ \begin{matrix}\longrightarrow\\[-3.5mm]\longrightarrow\\[-3.5mm]\longrightarrow\\[-3.5mm]\longrightarrow\\[-3.5mm]\longleftarrow\\[-3.5mm]
\longleftarrow\\[-3.5mm]\longleftarrow
\end{matrix}\ VP_3\ \begin{matrix}\longrightarrow\\[-3.5mm] \longrightarrow\\[-3.5mm]\longrightarrow\\[-3.5mm]
\longleftarrow\\[-3.5mm]\longleftarrow \end{matrix}\ VP_2\ \begin{matrix} \longrightarrow\\[-3.5mm]\longrightarrow\\[-3.5mm]
\longleftarrow \end{matrix}\ VP_1$$
with face and degeneracy homomorphisms
\begin{align*}
& d_i: \VAP_n=VP_{n+1}\to \VAP_{n-1}=VP_n,\ i=0,\dots n,\\
& s_i: \VAP_n=VP_{n+1}\to \VAP_{n+1}=VP_{n+2},\ i=0,\dots, n.
\end{align*}
Note that the vertex group $\VAP_0=VP_1$ is trivial. The first nontrivial group in $VP_*$ is $VP_2$, which is a free group of rank $2$ with standard generators $\lambda_{1,2}$ and $\lambda_{2,1}$. (The standard presentation of $VP_n$ will be reviewed in subsection~\ref{virt}.) There is a (unique) simplicial map
$$
f_{(\lambda_{1,2},\lambda_{2,1})}\colon S^1\vee S^1\longrightarrow \VAP_*
$$
such that $f_{(\lambda_{1,2},\lambda_{2,1})}((0,1)_1)=\lambda_{1,2}$ and $f_{(\lambda_{1,2},\lambda_{2,1})}((0,1)_2)=\lambda_{2,1}$, which induces a simplicial homomorphism
\begin{equation}\label{equation1.1}
\tilde{\Theta}_V\colon F[S^1\vee S^1]\longrightarrow \VAP_*.
\end{equation}
Here $(0,1)_1$ and $(0,1)_2$ are the non-degenerate $1$-simplices of the first copy and the second of $S^1$ in the wedge $S^1\vee S^1$, respectively. Let $T_*=\tilde{\Theta}_V(F[S^1\vee S^1])$ be the image of $\tilde{\Theta}_V$. Then $T_*$ is the smallest simplicial subgroup of $\VAP_*$ containing $\lambda_{1,2},\lambda_{2,1}$ with
$$
\tilde{\Theta}_V\colon F[S^1\vee S^1]_n\longrightarrow T_n\leq \VAP_n=VP_{n+1}
$$
an isomorphism for $n\leq 1$. Different from the case on classical braids, the simplicial homomorphism $\tilde{\Theta}_V$ is not injective with its first non-trivial kernel occurring in the case $n=2$. For determining the homotopy type of $T_*$, we need more detailed technical information. Observe that, for $n\geq 2$, $F[S^1\vee S^1]_n$ is the free group of rank $2n$ having a basis given by the non-base-point $n$-simplices of $S^1\vee S^1$ in the degenerate form of
$$
s_{n-1}s_{n-2}\cdots s_{k}\hat{s}_{k-1}s_{k-2}\cdots s_0(0,1)_1, \quad s_{n-1}s_{n-2}\cdots s_{k}\hat{s}_{k-1}s_{k-2}\cdots s_0(0,1)_2
$$
with $1\leq k\leq n$, where $\hat{s}_{k-1}$ means that the degenerate operation $s_{k-1}$ is omitted. Let
\begin{equation}\label{a_{kl}}
\begin{array}{rcl}
a_{k,n+1-k}&=&\tilde\Theta_V(s_{n-1}s_{n-2}\cdots s_{k}\hat{s}_{k-1}s_{k-2}\cdots s_0(0,1)_1)\\
&=&s_{n-1}s_{n-2}\cdots s_{k}\hat{s}_{k-1}s_{k-2}\cdots s_0\lambda_{1,2},\\
\end{array}
\end{equation}
\begin{equation}\label{b_{kl}}
\begin{array}{rcl}
b_{k,n+1-k}&=&\tilde\Theta_V(s_{n-1}s_{n-2}\cdots s_{k}\hat{s}_{k-1}s_{k-2}\cdots s_0(0,1)_2)\\
&=&s_{n-1}s_{n-2}\cdots s_{k}\hat{s}_{k-1}s_{k-2}\cdots s_0\lambda_{2,1}\\
\end{array}
\end{equation}
be the elements in $VP_{n+1}$ for $1\leq k\leq n$.

\begin{lem}~\cite{BMW1}\label{lemma1.1}
In group $VP_{n+1}$, the following commuting rules hold:
$$
a_{k,n+1-k}a_{l,n+1-l}=a_{l,n+1-l}a_{k,n+1-k} \textrm{ and } b_{k,n+1-k}b_{l,n+1-l}=b_{l,n+1-l}b_{k,n+1-k}
$$
for $1\leq k,l\leq n$.
\end{lem}

By this lemma, the simplicial homomorphism $\tilde \Theta_V$ factors through the free product $K(\Z,1)\ast K(\Z,1)$ as the simplicial quotient group of $F[S^1\vee S^1]$, where $K(\Z,1)=F[S^1]$ is the abelianization of $F[S^1]$. Let
\begin{equation}\label{equation1.2}
\Theta_V\colon K(\Z,1)\ast K(\Z,1)\rOnto T_*\leq \VAP_*
\end{equation}
be the resulting simplicial homomorphism. Unfortunately the simplicial homomorphism $\Theta_V$ is still not injective. However we can determine the homotopy type of the kernel of $\Theta_V$ by the following key lemma.

\begin{lem}~\cite{BMW1}\label{lemma1.2}
Let $K_*$ be the kernel of the simplicial epimorphism $\Theta_V\colon K(\Z,1)\ast K(\Z,1)\twoheadrightarrow T_*$. Then $K_*$ is contractible.
\end{lem}

As a consequence of this key lemma, the simplicial epimorphism $$\Theta_V\colon K(\Z,1)\ast K(\Z,1)\twoheadrightarrow T_*$$ is a homotopy equivalence from the long exact sequence of homotopy groups on fibrations with using the fact that any simplicial epimorphism between simplicial groups is a fibration. Recall that the classifying space of the free product $G\ast H$ of simplicial groups $G$ and $H$ is homotopy equivalent to the wedge of the classifying spaces of $G$ and $H$ by the Whitehead theorem ~\cite[Theorem 5]{Whitehead} with the formal statement given
in ~\cite[Proposition 4.3]{Kan-Thurston}.  Since the classifying space $K(\Z,1)$ is $K(\Z,2)\simeq \CP^\infty$, the geometric realization of $K(\Z,1)\ast K(\Z,1)$ is homotopy equivalent to $\Omega (\CP^\infty\vee \CP^\infty)$. From a generalization of the Hilton-Milnor Theorem given by B. Gray~\cite{Gray}, there is a homotopy decomposition
$$
\Omega (\CP^\infty\vee \CP^\infty)\simeq S^1\times S^1\times \Omega S^3.
$$
Before we give the following main result of the present paper, let us recall the definition of Moore homotopy group. Let $G_*$ be a simplicial group. The \textit{Moore cycles} $\mathrm{Z}_n(G_*)\leq G_n$ is defined by
$$
\mathrm{Z}_n(G_*)=\bigcap_{i=0}^n\mathrm{Ker}(d_i\colon G_n\to G_{n-1})
$$
and the \textit{Moore boundaries} $\mathcal{B}_n(G_*)\leq G_n$ is defined by
$$
\mathcal{B}_n(G_*)=d_0\left(\bigcap_{i=1}^{n+1}\mathrm{Ker}(d_i\colon G_{n+1}\to G_n)\right).
$$
Simplicial identities guarantees that $\mathcal{B}_n(G_*)$ is a (normal) subgroup of $\mathrm{Z}_n(G_*)$. The \textit{Moore homotopy group} $\pi_n(G_*)$ is defined by
$$
\pi_n(G_*)=\mathrm{Z}_n(G_*)/\mathcal{B}_n(G_*).
$$
It is a classical result due to J. C. Moore ~\cite{Moore} that $\pi_n(G_*)$ is isomorphic to the $n$-th homotopy group of the geometric realization of $G_*$.

\begin{thm}~\cite{BMW1}\label{theorem1.3}
The geometric realization of the simplicial group $T_*$ is homotopy equivalent to the loop space
$$
\Omega (\CP^\infty\vee \CP^\infty)\simeq S^1\times S^1\times \Omega S^3.
$$
In particular, the Moore homotopy group
$$\pi_n(T_*)\cong \pi_n(\Omega S^3)=\pi_{n+1}(S^3)$$ for $n\geq 2$.\hfill $\Box$
\end{thm}

Intuitively the subgroup $T_{n-1}\leq VP_n$ can be described as follows. An $n$-strand virtual pure braid $\beta$ is called \textit{cabled from $2$-strand braid} if there exists a $2$-strand virtual pure braid $\beta'$ such that
$$
\beta=s_{n-2}s_{n-3}\cdots s_{k}\hat{s}_{k-1}s_{k-2}\cdots s_0\beta'
$$
for some $1\leq k\leq n$. Namely $\beta$ is obtained from a cabling operation on $2$-strand virtual braid. An illustrative picture on the cabling on classical $2$-strand braids can be found in~\cite{CW}. The subgroup $T_{n-1}\leq VP_n$ consists of the $n$-strand virtual pure braids that can be written as a product of $n$-strand virtual pure braids cabled from $2$-strand braids.

The canonical inclusion $P_n\leq VP_n$ induces a simplicial homomorphism $\phi\colon AP_*\to VAP_*$. Let $\check T_*=\Theta(F[S^1])\leq AP_*$ be the image of $F[S^1]$ in $AP_*$. The restriction of $\phi$ to $\check T_*$ gives a simplicial homomorphism
$$
\phi|\colon \check T_*\longrightarrow T_*.
$$

\begin{thm}~\cite{BMW1}\label{theorem1.5}
The simplicial homomorphism $\phi|\colon \check T_*\longrightarrow T_*$ induces an isomorphism of homotopy groups
$$
\phi|_*\colon \pi_n(\check T_*)\rTo^{\cong} \pi_n(T_*)
$$
for $n\geq 2$.\hfill $\Box$
\end{thm}

\section{The cablings of virtual pure braid groups}\label{section3}

Our motivation on exploring the cabling presentation of virtual pure braid groups is to provide some fundamental information for proving Lemma~\ref{lemma1.2}, which is a key lemma for proving Theorem~\ref{theorem1.3}. The main point here is that, for handling the homotopy type of the simplicial group $K_*$ in Lemma~\ref{lemma1.2}, it is important to give a new presentation virtual pure braid groups in terms of the cabling virtual pure braids $a_{k,l}$ and $b_{s,t}$ given in formulae~(\ref{a_{kl}}) and ~(\ref{b_{kl}}).

\subsection{Virtual braid group} \label{virt}

The virtual braid group $VB_n$ is generated by elements
$$
\sigma_1,\,  \sigma_2, \, \ldots , \, \sigma_{n-1}, \, \rho_1, \, \rho_2, \, \ldots , \, \rho_{n-1},
$$
where $\sigma_1,\,  \sigma_2, \, \ldots , \, \sigma_{n-1}$ generate the classical braid group $B_n$ and
the elements $\rho_1$,  $\rho_2$,  $\ldots $,  $\rho_{n-1}$ generate the symmetric group
$S_n$. Hence, $VB_n$ is defined by relations of $B_n$, relations of $S_n$
and mixed relation:
$$
\sigma_i \rho_j = \rho_j \sigma_i,~~~|i-j| > 1,
$$
$$
\rho_i \rho_{i+1} \sigma_i = \sigma_{i+1} \rho_i \rho_{i+1}.
$$

As for the classical braid groups there exists the canonical
epimorphism of $VB_n$ onto the symmetric group $VB_n\twoheadrightarrow S_n$ with the
kernel called the \textit{virtual pure  braid group} $VP_n$. So we have a
short exact sequence
\begin{equation*}
1 \to VP_n \to VB_n \to S_n \to 1.
\end{equation*}
Define the following elements in $VP_n$:
$$
\lambda_{i,i+1} = \rho_i \, \sigma_i^{-1},~~~
\lambda_{i+1,i} = \rho_i \, \lambda_{i,i+1} \, \rho_i = \sigma_i^{-1} \, \rho_i,
~~~i=1, 2, \ldots, n-1,
$$
$$
\lambda_{ij} = \rho_{j-1} \, \rho_{j-2} \ldots \rho_{i+1} \, \lambda_{i,i+1} \, \rho_{i+1}
\ldots \rho_{j-2} \, \rho_{j-1},
$$
$$
\lambda_{ji} = \rho_{j-1} \, \rho_{j-2} \ldots \rho_{i+1} \, \lambda_{i+1,i} \, \rho_{i+1}
\ldots \rho_{j-2} \, \rho_{j-1}, ~~~1 \leq i < j-1 \leq n-1.
$$
It is shown in \cite{B} that the group $VP_n\ (n\geq 2)$ admits a
presentation with the  generators $\lambda_{ij},\ 1\leq i\neq j\leq n,$
and the following relations:
\begin{align}
& \lambda_{ij}\lambda_{kl}=\lambda_{kl}\lambda_{ij} \label{rel},\\
&
\lambda_{ki}\lambda_{kj}\lambda_{ij}=\lambda_{ij}\lambda_{kj}\lambda_{ki}
\label{relation},
\end{align}
where distinct letters stand for distinct indices.

\subsection{The cablings of virtual pure braid groups}

By using the same ideas in the work~\cite{BCWW,CW} on the classical braids, we have a simplcial group
$$
\VAP_* :\ \ \ \ldots\ \begin{matrix}\longrightarrow\\[-3.5mm] \ldots\\[-2.5mm]\longrightarrow\\[-3.5mm]
\longleftarrow\\[-3.5mm]\ldots\\[-2.5mm]\longleftarrow \end{matrix}\ VP_4 \ \begin{matrix}\longrightarrow\\[-3.5mm]\longrightarrow\\[-3.5mm]\longrightarrow\\[-3.5mm]\longrightarrow\\[-3.5mm]\longleftarrow\\[-3.5mm]
\longleftarrow\\[-3.5mm]\longleftarrow
\end{matrix}\ VP_3\ \begin{matrix}\longrightarrow\\[-3.5mm] \longrightarrow\\[-3.5mm]\longrightarrow\\[-3.5mm]
\longleftarrow\\[-3.5mm]\longleftarrow \end{matrix}\ VP_2\ \begin{matrix} \longrightarrow\\[-3.5mm]\longrightarrow\\[-3.5mm]
\longleftarrow \end{matrix}\ VP_1$$
on pure virtual braid groups with $\VAP_n=VP_{n+1}$, the face homomorphism
$$
d_i : \VAP_n=VP_{n+1} \longrightarrow \VAP_{n-1}=VP_n
$$
given by deleting $(i+1)$th strand for $0\leq i\leq n$, and the degeneracy homomorphism
$$
s_i : \VAP_n=VP_{n+1} \longrightarrow \VAP_{n+1}=VP_{n+2}
$$
given by doubling the $(i+1)$th strand for $0\leq i\leq n$.

The following proposition is obtained from the geometric description, which can be regarded as the formal definition. The proof of the following proposition is straightforward.

\begin{prop} \label{p3.1}
The sequence of groups $\VAP_*$ with $\VAP_n=VP_{n+1}$ for $n\geq0$ is a simplicial group under the faces
$d_i : \VAP_{n-1}=VP_{n} \longrightarrow \VAP_{n-2}=VP_{n-1}$, $0\leq i\leq n-1$, and
 degeneracies $s_i : \VAP_{n-1}=VP_{n} \longrightarrow \VAP_{n}=VP_{n+1}$, $0\leq i\leq n-1$, given the group homomorphism with acting on the generators $\lambda_{k,l}$ and $\lambda_{l,k}$, $1 \leq k < l \leq n$, of $VP_{n}$ by the rules
$$
s_i (\lambda_{k,l}) = \left\{
\begin{array}{lcl}
\lambda_{k+1,l+1} & \textrm{ for }&i < k-1,\\
\lambda_{k,l+1} \lambda_{k+1,l+1} & \textrm{ for } &i = k-1, \\
\lambda_{k,l+1} & \textrm{ for } &k-1 < i < l-1,\\
& \\
\lambda_{k,l+1}^{\lambda_{1l} \lambda_{2l} \ldots \lambda_{k-1,l}} \, \lambda_{k,l} & \textrm{ for }& i = l-1, \\
& \\
\lambda_{k,l} & \textrm{ for } &i > l-1,
\end{array}
\right.
$$
$$
s_i (\lambda_{l,k}) = \left\{
\begin{array}{lcl}
\lambda_{l+1,k+1} &\textrm{ for }&i < k-1,\\
\lambda_{l+1,k+1} \lambda_{l+1,k} &\textrm{ for }&i = k-1, \\
\lambda_{l+1,k} & \textrm{ for }&k-1 < i < l-1,\\
& \\
\lambda_{l,k} \, \lambda_{l+1,k}^{\lambda_{l1}^{-1} \lambda_{l2}^{-1} \ldots \lambda_{l,k-1}^{-1}}   & \textrm{ for }&i = l-1, \\
& \\
\lambda_{l,k} & \textrm{ for }&i > l-1,
\end{array}
\right.
$$
$$
d_i (\lambda_{k,l}) = \left\{
\begin{array}{lcl}
\lambda_{k-1,l-1} & \textrm{ for } &0 \leq i < k-1, \\
1 & \textrm{ for } & i =k-1,\\
\lambda_{k,l-1} & \textrm{ for } &k-1 < i < l-1, \\
1 & \textrm{ for } & i =l,\\
\lambda_{k,l} & \textrm{ for } &l-1 < i \leq n-1, \\
\end{array}
\right.
$$
$$
d_i (\lambda_{l,k}) = \left\{
\begin{array}{lcl}
\lambda_{k-1,l-1} & \textrm{ for } &0 \leq i < k-1, \\
1 & \textrm{ for } & i = k-1,\\
\lambda_{l-1,k} & \textrm{ for } &k-1 < i < l-1, \\
1 & \textrm{ for } & i = l-1,\\
\lambda_{l,k} & \textrm{ for } &l-1 < i \leq n-1, \\
\end{array}
\right.
$$
where $y^x=x^{-1}yx$. \hfill $\Box$
\end{prop}

The simplicial subgroup $T_*$ has been defined in Section~\ref{section2} as the image of $\tilde\Theta_V$ in equation~(\ref{equation1.1}). It is routine to see that the group $T_n$ as a subgroup of $VP_{n+1}$ can be constructed recursively as follows:
\begin{enumerate}
\item[] $T_0=\{1\}$, $T_1=VP_2$, and $T_{n+1} = \langle s_0(T_n), s_1(T_n), \ldots,
s_{n}(T_n) \rangle. $
\end{enumerate}

Let $a_{k,l}$ and $b_{s,t}$ be defined in equations~(\ref{a_{kl}}) and ~(\ref{b_{kl}}), respectively. We can prove the following formulae.
\begin{equation}\label{equation3.1}
a_{n-k,k} = \left\{
\begin{array}{lr}
\lambda_{1n} \lambda_{2n} \ldots \lambda_{n-1,n} & for  ~k = 1,\\
\lambda_{1n} \lambda_{2n} \ldots \lambda_{n-k,n} a_{n-k,k-1} & for ~1 < k < n, \\
\lambda_{1n} a_{1,n-2} & for  ~k = n-1,
\end{array}
\right.
\end{equation}
\begin{equation}\label{equation3.2}
b_{n-k,k} = \left\{
\begin{array}{lr}
\lambda_{n,n-1} \lambda_{n,n-2} \ldots \lambda_{n1} & for  ~k = 1,\\
 b_{n-k,k-1} \lambda_{n,n-k} \lambda_{n,n-k-1} \ldots \lambda_{n1} & for ~1 < k < n, \\
b_{1,n-2}  \lambda_{n1} & for  ~k = n-1,
\end{array}
\right.
\end{equation}
Moreover the generators $\lambda_{ij}$ can be written in terms of $a_{k,l}$ and $b_{s,t}$ as follows:
\begin{equation}\label{equation3.3}
\lambda_{kn} = \left\{
\begin{array}{lr}
a_{1,n-1} \, a_{1,n-2}^{-1} & for  ~k = 1,\\
a_{k-1,n-k} \, a_{k-1,n-k+1}^{-1} \, a_{k,n-k} \, a_{k,n-k-1}^{-1} & for ~1 < k < n, \\
a_{n-2,1} \, a_{n-2,2}^{-1} \, a_{n-1,1} & for  ~k = n-1,
\end{array}
\right.
\end{equation}
\begin{equation}\label{equation3.4}
\lambda_{nk} = \left\{
\begin{array}{lr}
b_{1,n-2}^{-1} \, b_{1,n-1} & for  ~k = 1,\\
b_{k,n-k-1}^{-1} \, b_{k,n-k} \, b_{k-1,n-k+1}^{-1} \, b_{k-1,n-k} & for ~1 < k < n, \\
b_{n-1,1} \, b_{n-2,2}^{-1} \, b_{n-2,1} & for  ~k = n-1,
\end{array}
\right.
\end{equation}

From the above formulae, we have the following proposition.

\begin{prop} \label{p3.3}
Consider $VP_k$ as a subgroup of $VP_{k+1}$ by adding a trivial strand in the end. Then
\begin{enumerate}
\item  The subgroup $T_{n-1}$ of $VP_n$, $n \geq 3$, is generated by elements $a_{kl}$, $b_{kl}$, $k+l = n$.
\item The group $VP_n=\la T_1, T_2,\ldots, T_n\ra$ generated by $a_{k,l}$ and $b_{k,l}$ for $2\leq k+l\leq n, 1\leq k,l\leq n-1$.
\item $VP_{n+1}=\la VP_n, s_0VP_n, s_1 VP_n,\ldots, s_{n-1}VP_n \ra$ for $n\geq 2$.
\end{enumerate}
\end{prop}

By this proposition $VP_n$ admits a new generating system in terms of $a_{k,l}$ and $b_{k,l}$, which are intuitively given by cabling on $2$-strand braid and adding trivial strands to the end. Moreover, we have a presentation of $VP_n$ in terms of $a_{k,l}$ and $b_{s,t}$ from formulae ~(\ref{equation3.3}) and~(\ref{equation3.4}).

Let $n\geq 4$. Let $\mathcal{R}^V(n)$ denote the defining relations~(\ref{rel}) and ~(\ref{relation}) of $VP_n$. By applying the doubling homomorphism $s_t\colon VP_n\to VP_{n+1}$ to $\mathcal{R}^V(n)$ , we have the following equations
\begin{align}
& s_t(\lambda_{ij})s_t(\lambda_{kl})=s_t(\lambda_{kl})s_t(\lambda_{ij}) \label{equation3.7},\\
&
s_t(\lambda_{ki})s_t(\lambda_{kj})s_t(\lambda_{ij})=s_t(\lambda_{ij})s_t(\lambda_{kj})s_t(\lambda_{ki})
\label{equation3.8}
\end{align}
in $VP_{n+1}$ for $1\leq i,j,k,l\leq n$ with distinct letters standing for distinct indices, which is denoted as $s_t(\mathcal{R}^V(n))$.

\begin{thm}~\cite{BMW2}
Let $n\geq 4$. Consider $VP_n$ as a subgroup of $VP_{n+1}$ by adding a trivial strand in the end. Then
$$
\mathcal{R}^V(n)\cup\bigcup_{i=0}^{n-1}s_i(\mathcal{R}^V(n))
$$
gives the full set of the defining relations for $VP_{n+1}$.\hfill $\Box$
%
%
%
%
%
%
%
%
%
%
%
%
%
%
%
%
%
%
\end{thm}

A significance of this theorem is that the group structure of $VP_n$ with $n\geq 5$ is determined by $VP_3$, $VP_4$ and virtual cablings given by iterated degeneracy operations $s_{j_t}s_{j_{t-1}}\cdots s_{j_1}$ on the generators and defining relations.

\section{Applications to the theory of virtual Brunnian braids}

In the theory of classical braids, a \textit{Brunnian braid} means a pure braid that becomes trivial after
removing any one of its strands. A typical example of a
$3$-strand Brunnian braid on a disk is the braid given by the
expression $(\sigma_1^{-1}\sigma_2)^3$, where $\sigma_1$ and
$\sigma_2$ are the standard generators of the 3-strand braid group
$\langle \sigma_1,\sigma_2\ |\
\sigma_1\sigma_2\sigma_1=\sigma_2\sigma_1\sigma_2\rangle$. A classical question proposed by G.~S.~Makanin~\cite{Makanin} in 1980 is to determine a set of
generators for Brunnian braids over the disk. Brunnian braids were called
\textit{smooth braids} by Makanin. This question was
answered by D.~L.~Johnson~\cite{Johnson} and G.~G.~Gurzo~\cite{Gurzo}.
A different approach to
this question can be found in~\cite{LW, Wu1}.
In the 1970s, H.~W.~Levinson \cite{Lev1, Lev2} defined a notion of {\it
$k$-decomposable} braid. It means a braid which becomes  trivial
after removal of any arbitrary $k$ strings. In his terminology a
{\it decomposable} braid means {\it $1$-decomposable} and therefore, Brunnian.

A connection between Brunnian braids and the homotopy groups of spheres
was given in~\cite{BCWW}. In particular, the following exact sequence
\begin{equation}
1\to \Brun_{n+1}(S^2)\to \Brun_n(D^2) \to \Brun_n(S^2)\to \pi_{n-1}(S^2)
\to 1
\label{eq:S2D2}
\end{equation}
was proved for $n>4$.

 In her
book~\cite[Question 23, p. 219]{Birman}, Birman asked how to determine a free basis for $\Brun_n(D^2)\cap R_{n-1}$ where
$$R_{n-1}=\Ker(B_n(D^2)\to B_n(S^2)).$$
Her motivation was that the kernel of the Gassner representation
is a subgroup of $\Brun_n(D^2)\cap R_{n-1}$.
From the exact sequence (\ref{eq:S2D2}) it follows that Birman's
question, for $n>5$, is about a free basis of
Brunnian braids over the sphere $S^2$.
As far as we know this  question remains open.

In the theory of virtual braids, we call an $n$-strand virtual pure braid \textit{Brunnian} if it becomes trivial after removing any one of its strands. Let $\Brun^V_n$ denote the group of $n$-strand virtual Brunnian braids.  Let $\Brun_n\leq P_n$ be the group of $n$-strand classical Brunnian braids. For a subgroup $H\leq VP_n$, we call a braid $\beta$ is \textit{virtual Brunnian in $H$} if $\beta\in H\cap \Brun^V_n$, and \textit{classical Brunnian in $H$} if $\beta\in H\cap \Brun_n$.

A natural and fundamental question in the theory of Brunnian braids is how to determine virtual or classical Brunnian braids in a given subgroup $H$. This question is a generalization of Birman's question~\cite[Question 23, p. 219]{Birman}. The constraint subgroups $H$ considered in this article are $T_{n-1}$ and $\check T_{n-1}$, namely the subgroups generated by $n$-strand braids cabled from $2$-strand virtual and classical braids, respectively. Theorems~\ref{theorem1.3} and~\ref{theorem1.5} admit applications in these cases.

The Moore boundaries $\mathcal{B}_{n-1}(T_*)$ is understandable in the following sense. Since $T_*$ is a simplicial quotient group of $F[S^1\vee S^1]$, the simplicial epimorphism $F[S^1\vee S^1]\to T_*$ induces an epimorphism
$$
\mathcal{B}_{n-1}(F[S^1\vee S^1])\twoheadrightarrow \mathcal{B}_{n-1}(T_*).
$$
An explicit construction of a generating set for $\mathcal{B}_{n-1}(F[S^1\vee S^1])$ given in~\cite{Wu1} gives a generating set for $\mathcal{B}_{n-1}(T_*)$. More precisely, let
\begin{equation}\label{equation4.2}
\begin{array}{cccccc}
y_1=a_{1,n-1}^{-1}& y_2=a_{1,n-1}a_{2,n-2}^{-1} & y_3=a_{2,n-2}a_{3,n-3}^{-1}&\cdots& y_{n-1}=a_{n-2,2}a_{n-1,1}^{-1}& y_{n}=a_{n-1,1},\\
z_1=b_{1,n-1}^{-1}& z_2=b_{1,n-1}b_{2,n-2}^{-1} & z_3=b_{2,n-2}b_{3,n-3}^{-1}&\cdots& z_{n-1}=b_{n-2,2}b_{n-1,1}^{-1}& z_{n}=b_{n-1,1}\\
\end{array}
\end{equation}
be elements in $T_{n-1}\leq VP_n$. Let $R_i=\la\la y_i,z_i\ra\ra^{T_{n-1}}$ be the normal closure of $y_i$ and $z_i$ in $T_{n-1}$ for $1\leq i\leq n$. Let
$$
[R_1,R_2,\ldots,R_n]_S=\prod_{\sigma\in\Sigma_n}[[[R_{\sigma(1)}, R_{\sigma(2)}], R_{\sigma(3)}], \ldots, R_{\sigma(n)}]
$$
be the symmetric commutator subgroup of the subgroups $R_1, R_2,\ldots,R_n$.
Then, by~\cite[Theorem 4.4]{Wu1} and~\cite[Theorem 1.1]{LW}, the Moore boundaries
\begin{equation}\label{equation4.3}
\mathcal{B}_{n-1}(T_*)=[R_1,R_2,\ldots,R_n]_S.
\end{equation}

Theorem ~\ref{theorem1.3} can be reformulated as
$$
[R_1,R_2,\ldots,R_n]_S\unlhd T_{n-1}\cap \Brun^V_n
$$
with cokernel given by $\pi_n(S^3)$ for $n\geq 3$. The braided interpretation is that the homotopy group $\pi_n(S^3)$ measures the virtual Brunnian braids in $T_{n-1}$ away from the symmetric commutator subgroup $[R_1,R_2,\ldots,R_n]_S$. This interpretation extends the main result in ~\cite{CW} on classical Brunnian braids to virtual Brunnian braids.

The Moore boundaries $\mathcal{B}_{n-1}(\check T_*)$ can be described as a symmetric commutator subgroup by the same reasons. For describing it as a subgroup of $VP_n$, let
$$
c_{1,1}=a_{1,1}^{-1}b_{1,1}^{-1}=\sigma_1\rho_1\rho_1\sigma_1=\sigma_1^2
$$
and
\begin{equation}\label{c_{kl}}
c_{k,n-k}=s_{n-2}s_{n-3}\cdots s_{k}\hat{s}_{k-1}s_{k-2}\cdots s_0c_{1,1}
\end{equation}
for $1\leq k\leq n-1$ with $n\geq 3$. Let
\begin{equation}\label{equation4.5}
w_k=c_{k-1,n-k+1}c_{k,n-k}^{-1}
\end{equation}
for $1\leq k\leq n$, where $c_{0,n}=c_{n,0}=1$. Let $\check R_k=\la\la w_k\ra\ra^{\check T_{n-1}}$ be the normal closure of $w_k$ in $\check T_{n-1}$. Then, by~\cite[Theorem 4.4]{Wu1} and~\cite[Theorem 1.1]{LW}, the Moore boundaries
\begin{equation}\label{equation4.6}
\mathcal{B}_{n-1}(\check T_*)=[\check R_1,\check R_2,\ldots,\check R_n]_S.
\end{equation}

Theorem ~\ref{theorem1.3} together with Theorem~\ref{theorem1.5}, we have the following theorem in the theory of Brunnian braids.

\begin{thm}~\cite{BMW1}\label{theorem1.6}
\begin{enumerate}
\item \textbf{Exponent Property.} Let $\beta\in T_{n-1}\cap \Brun^V_n$ be a virtual Brunnian braid in $T_{n-1}$ with $n\geq 4$. Then  there exists an odd integer $k$ given as a product of distinct prime integers such that
$$
\beta^{4k}\in [R_1,R_2,\ldots,R_n]_S.
$$
\item \textbf{Non-triviality Property.} $[R_1,R_2,\ldots,R_n]_S$ is a proper normal subgroup of $T_{n-1}\cap \Brun^V_n$ for $n\geq 2$.
\item \textbf{Decomposition Property.} Let $\beta\in T_{n-1}\cap \Brun^V_n$ be a virtual Brunnian braid in $T_{n-1}$ with $n\geq 3$. Then there exists a classical Brunnian braid $\tilde \beta$ in $\check T_{n-1}$ (that is $\tilde \beta\in \check T_{n-1}\cap \Brun_n$) and $\gamma\in [R_1,R_2,\ldots,R_n]_S$ such that
$$
\beta=\tilde\beta \cdot \gamma.
$$
Moreover, if $$\tilde\beta_1\gamma_1=\tilde\beta_2\gamma_2$$ with $\tilde\beta_1,\tilde\beta_2\in \check T_{n-1}\cap \Brun_n$ and $\gamma_1,\gamma_2\in [R_1,R_2,\ldots,R_n]_S$, then there exists $\delta\in [\check R_1,\check R_2,\ldots,\check R_n]_S$ such that
$$\tilde \beta_2=\tilde\beta_1\cdot \delta\textrm{ and }\gamma_2=\delta^{-1}\gamma_1.$$\hfill $\Box$
\end{enumerate}
\end{thm}

The exponent and non-triviality properties are interpretations  of the deep theorems on the exponents and non-triviality of the homotopy groups of $S^3$ into the content of virtual braids. It seems hard to have a direct proof of the exponent property using current techniques in the theory of braids. The geometric exploration on the exponent property of virtual Brunnian braids might give new insight in the topic.

The decomposition property follows from analysing the homotopy behalves of the simplicial embedding $AP_*\to \VAP_*$. A significance of the decomposition property is that the study on virtual Brunnian braids in $T_{n-1}$ can be reduced to the study on classical Brunnian braids in $\check T_{n-1}$ and the symmetric commutator subgroup $[R_1,R_2,\ldots,R_n]_S$.

\end{document}